\documentclass{birkjour}
\usepackage{amsmath, amsthm, amscd, amsfonts, amssymb, graphicx, color}
 \newtheorem{thm}{Theorem}[section]
 \newtheorem{cor}[thm]{Corollary}
  \newtheorem{con}[thm]{Conjecture}
 \newtheorem{lem}[thm]{Lemma}
 \newtheorem{prop}[thm]{Proposition}
 \theoremstyle{definition}
 \newtheorem{defn}[thm]{Definition}
 \theoremstyle{remark}
 \newtheorem{rem}[thm]{Remark}
 
 \numberwithin{equation}{section}

\begin{document}

%
%
%
%
%
%
%
%
%

\title[{Langlands reciprocity  for $C^*$-algebras}]
 {Langlands reciprocity for $C^*$-algebras}

\author[Nikolaev]{Igor ~V. ~Nikolaev}

\address{%
Department of Mathematics \\
St.~John's University\\
 8000 Utopia Parkway\\
New York,  NY 11439\\
United States}

\email{igor.v.nikolaev@gmail.com}




\date{October 11, 2019}

\begin{abstract}
We introduce  a $C^*$-algebra ${\mathcal A}_V$ of  a variety $V$ over the number field $K$
and a  $C^*$-algebra ${\mathcal A}_G$ of a reductive group $G$
over the ring of adeles of $K$.   Using Pimsner's Theorem,  we construct 
 an embedding ${\mathcal A}_V\hookrightarrow {\mathcal A}_G$,
 where  $V$ is a $G$-coherent variety, e.g. the Shimura variety of $G$. 
The embedding is an  analog of the Langlands reciprocity 
for $C^*$-algebras.
It follows from  the $K$-theory of  the  inclusion ${\mathcal A}_V\subset{\mathcal A}_G$
 that  the Hasse-Weil  $L$-function  of  $V$ 
is  a product of the automorphic $L$-functions
corresponding  to irreducible representations of the group $G$. 
\end{abstract}

\maketitle

\section{Introduction}
The Langlands conjectures say   that all zeta functions are automorphic
[Langlands 1978]  \cite{Lan1}.  In this note we  study   (one of)  the conjectures 
in terms of  the  $C^*$-algebras [Dixmier 1977]  \cite{D}. 
Namely,  denote by $G({\mathbf A}_K)$ a reductive group $G$  over the ring of
adeles  ${\mathbf A}_K$ of a number field $K$ and by $G(K)$ 
its discrete subgroup over $K$. 
 The  Banach algebra $L^1(G(K) \backslash G({\mathbf A}_K))$ consists of 
the  integrable complex-valued  functions endowed  with the operator norm.  The addition of functions 
$f_1, f_2\in L^1(G(K) \backslash G({\mathbf A}_K))$ is defined pointwise  and 
 multiplication is given by the convolution product:
\begin{equation}\label{eq1}
(f_1\ast f_2)(g)=\int_{G(K)\backslash G({\mathbf A}_{K})} f_1(gh^{-1})f_2(h)dh. 
\end{equation}  
Consider the enveloping $C^*$-algebra,   ${\mathcal A}_G$,  of  the algebra  $L^1(G(K) \backslash G({\mathbf A}_K) )$;
we refer the reader to  [Dixmier 1977,  Section 13.9]  \cite{D}
for details of this construction.    The algebra ${\mathcal A}_G$ encodes all  unitary   irreducible    representations of 
the locally compact group $G({\mathbf A}_K)$  induced by   $G(K)$.
Such representations are related to   the automorphic cusp forms and non-abelian class 
field theory [Gelbart 1984] \cite{Gel1}.  The algebra  ${\mathcal A}_G$ has an amazingly
simple structure.  Namely,  let us assume  $G\cong GL_n$. 
Then ${\mathcal A}_G$  is a  stationary {\it AF-algebra},  see lemma \ref{lem3.1}; 
such an algebra is defined  by a positive integer matrix $B\in SL_n({\bf Z})$
[Bratteli  1972]    \cite{Bra1}  and its $K$-theory is  well understood  [Effros 1981]  \cite{E}.

Let $V$ be a complex projective variety.
For an automorphism $\sigma: V\to V$ and an invertible
sheaf $\mathcal{L}$ of the linear forms on $V$,  one can construct  
a twisted homogeneous coordinate ring $B(V, \mathcal{L}, \sigma)$ of the variety $V$, 
i.e. a non-commutative ring such that:
\begin{equation}
Mod~(B(V, \mathcal{L}, \sigma)) ~/~Tors \cong ~Coh~(V),  
\end{equation}
where $Mod$ is the category of graded left modules over the graded ring 
$B(V, \mathcal{L}, \sigma)$, $Tors$ the full subcategory of $Mod$ of the 
torsion modules and $Coh$ the category of quasi-coherent sheaves on the 
variety  $V$ [Stafford \& van ~den ~Bergh 2001]  \cite[p.180]{StaVdb1}.    
The norm-closure of a self-adjoint representation
of the ring $B(V, \mathcal{L}, \sigma)$ by  the  linear operators on 
a Hilbert space  is  called the {\it Serre $C^*$-algebra}   of 
$V$ \cite{Nik1}.  In what follows,  we shall  focus on the case 
when $V$ is defined over a number field $K$, i.e. $V$ is an arithmetic variety.
  The corresponding  Serre $C^*$-algebra is denoted by  ${\mathcal A}_V$.  The Hasse-Weil 
$L$-function of $V$ was calculated in \cite{Nik1}  in terms of the $K$-theory of algebra  ${\mathcal A}_V$.

It is known that  the  Langlands philosophy  does not distinguish
between  the arithmetic  and automorphic objects [Langlands 1978]  \cite{Lan1}. 
Therefore one can expect a regular map between  the $C^*$-algebras ${\mathcal A}_V$ and  ${\mathcal A}_G$,
provided $V$ is a $G$-coherent variety, see Definition \ref{dfn0}. 
We prove that  such a map is an embedding  ${\mathcal A}_V\hookrightarrow  {\mathcal A}_G$.
To give an exact statement,  we shall need  the following notions. 
The   $i$-th  {\it trace cohomology} $\{H^i_{tr}(V)~|~0\le i\le 2\dim_{\mathbf{C}}V\}$
of an  arithmetic   variety $V$ is an additive  abelian subgroup of $\mathbf{R}$ obtained 
from a canonical trace on the Serre $C^*$-algebra of $V$ \cite{Nik1}. 
Likewise,  the group $K_0({\mathcal A}_G)$ of  the stationary AF-algebra ${\mathcal A}_G$ 
is  an additive abelian subgroup  of  $\mathbf{R}$  [Effros 1981]  \cite[Chapter 6]{E}. 
\begin{defn}\label{dfn0}
The arithmetic variety $V$ is called $G$-coherent, if
\begin{equation}\label{eq4}
H^i_{tr}(V)\subseteq K_0({\mathcal A}_G) \quad\hbox{for all}  \quad 0\le i\le 2\dim_{\mathbf{C}}V.
\end{equation}
\end{defn}
\begin{rem}
If $V\cong Sh~(G,X)$ is the Shimura variety corresponding to  the Shimura datum $(G, X)$ [Deligne 1971]  \cite{Del1}, 
then $V$ is a $G$-coherent  variety.  This remark  follows from an adaption of the argument for the Shimura curves
considered in  Section 4.    To put it simple, the arithmetic variety $V$ is $G$-coherent if
for all  $0\le i\le 2\dim_{\mathbf{C}}V$  the number fields $\mathbf{k}_i:=H^i_{tr}(V)\otimes\mathbf{Q}$ are 
subfields of (or coincide with) a number field  $\mathbf{K}:=K_0({\mathcal A}_G)\otimes\mathbf{Q}$.
A quick example are elliptic curves with complex multiplication, see Proposition \ref{prp2}. 
\end{rem}
\begin{thm}\label{thm2}
There exists a canonical embedding  ${\mathcal A}_V\hookrightarrow  {\mathcal A}_G$,
where $V$ is a $G$-coherent variety. 
\end{thm}
\begin{rem}
Theorem \ref{thm2} can be viewed as  an analog of the  Langlands 
 reciprocity for $C^*$-algebras.  In other words,   the coordinate ring  ${\mathcal A}_V$
 of a $G$-coherent   variety  $V$ is a sub-algebra of  the algebra  ${\mathcal A}_G$.  
 \end{rem}
An application of theorem \ref{thm2} is as follows. 
Recall that to each arithmetic variety $V$ one can attach the Hasse-Weil  (motivic) 
$L$-function.  Likewise,  to each irreducible representation  of the
 group $G({\mathbf A}_{K})$ one can attach an automorphic (standard)  $L$-function, see  
 [Gelbart 1984] \cite{Gel1} and [Langlands 1978] \cite{Lan1}.
 Theorem \ref{thm2} implies one of the conjectures of 
[Langlands 1978]  \cite{Lan1}. 
\begin{cor}\label{cor1}
The Hasse-Weil   $L$-function of  a $G$-coherent  variety $V$
is a product of  the automorphic  $L$-functions. 
\end{cor}
The paper is organized as follows.  The definitions and preliminary results can be found in Section 2. 
Theorem  \ref{thm2}  and corollary \ref{cor1} are  proved in Section 3.
An example is  constructed   in Section 4.

\section{Preliminaries}
This section is a brief account of preliminary facts involved in our paper;
we refer the reader to  [Bratteli  1972]  \cite{Bra1},  [Dixmier 1977]  \cite{D}, 
[Langlands 1978] \cite{Lan1} and  [Stafford \& van ~den ~Bergh 2001]  \cite{StaVdb1}.

\subsection{AF-algebras}
A {\it $C^*$-algebra} is an algebra $A$ over $\mathbf{C}$ with a norm
$a\mapsto ||a||$ and an involution $a\mapsto a^*$ such that
it is complete with respect to the norm and $||ab||\le ||a||~ ||b||$
and $||a^*a||=||a^2||$ for all $a,b\in A$.
Any commutative $C^*$-algebra is  isomorphic
to the algebra $C_0(X)$ of continuous complex-valued
functions on some locally compact Hausdorff space $X$; 
otherwise, $A$ represents a noncommutative  topological
space.

An {\it AF-algebra}  (Approximately Finite $C^*$-algebra) is defined to
be the  norm closure of an ascending sequence of   finite dimensional
$C^*$-algebras $M_n$,  where  $M_n$ is the $C^*$-algebra of the $n\times n$ matrices
with entries in $\mathbf{C}$. Here the index $n=(n_1,\dots,n_k)$ represents
the  semi-simple matrix algebra $M_n=M_{n_1}\oplus\dots\oplus M_{n_k}$.
The ascending sequence mentioned above  can be written as 
\begin{equation}\label{eq2}
M_1\buildrel\rm\varphi_1\over\longrightarrow M_2
   \buildrel\rm\varphi_2\over\longrightarrow\dots,
\end{equation}
where $M_i$ are the finite dimensional $C^*$-algebras and
$\varphi_i$ the homomorphisms between such algebras.  
If $\varphi_i=Const$, then the AF-algebra ${\mathcal A}$ is called 
{\it stationary};  such an algebra defines and is defined by 
a {\it shift automorphism}  $\sigma_{\varphi}: {\mathcal A}\to {\mathcal A}$
corresponding to a map $i\mapsto i+1$ on  $\varphi_i$
[Effros 1981]  \cite[p.37]{E}. 
The homomorphisms $\varphi_i$ can be arranged into  a graph as follows. 
Let  $M_i=M_{i_1}\oplus\dots\oplus M_{i_k}$ and 
$M_{i'}=M_{i_1'}\oplus\dots\oplus M_{i_k'}$ be 
the semi-simple $C^*$-algebras and $\varphi_i: M_i\to M_{i'}$ the  homomorphism. 
One has  two sets of vertices $V_{i_1},\dots, V_{i_k}$ and $V_{i_1'},\dots, V_{i_k'}$
joined by  $b_{rs}$ edges  whenever the summand $M_{i_r}$ contains $b_{rs}$
copies of the summand $M_{i_s'}$ under the embedding $\varphi_i$. 
As $i$ varies, one obtains an infinite graph called the  {\it Bratteli diagram} of the
AF-algebra.  The matrix $B=(b_{rs})$ is known as  a {\it partial multiplicity} matrix;
an infinite sequence of $B_i$ defines a unique AF-algebra.

For a unital $C^*$-algebra $A$, let $V(A)$
be the union (over $n$) of projections in the $n\times n$
matrix $C^*$-algebra with entries in $A$;
projections $p,q\in V(A)$ are {\it equivalent} if there exists a partial
isometry $u$ such that $p=u^*u$ and $q=uu^*$. The equivalence
class of projection $p$ is denoted by $[p]$;
the equivalence classes of orthogonal projections can be made to
a semigroup by putting $[p]+[q]=[p+q]$. The Grothendieck
completion of this semigroup to an abelian group is called
the  $K_0$-group of the algebra $A$.
The functor $A\to K_0(A)$ maps the category of unital
$C^*$-algebras into the category of abelian groups, so that
projections in the algebra $A$ correspond to a positive
cone  $K_0^+\subset K_0(A)$ and the unit element $1\in A$
corresponds to an order unit $u\in K_0(A)$.
The ordered abelian group $(K_0,K_0^+,u)$ with an order
unit  is called a {\it dimension group};  an order-isomorphism
class of the latter we denote by $(G,G^+)$.

If ${\mathcal A}$ is an AF-algebra, then its dimension group
$(K_0({\mathcal A}), K_0^+({\mathcal A}), u)$ is a complete isomorphism
invariant of algebra ${\mathcal A}$ [Elliott 1976]   \cite{Ell1}.  
The order-isomorphism  class $(K_0({\mathcal A}), K_0^+({\mathcal A}))$
 is an invariant of the {\it Morita equivalence} of algebra 
${\mathcal A}$,  i.e.  an isomorphism class in the category of 
finitely generated projective modules over ${\mathcal A}$.

\subsection{Trace cohomology}
Let $V$ be an $n$-dimensional complex  projective variety endowed with an automorphism $\sigma:V\to V$  
 and denote by $B(V, \mathcal{L}, \sigma)$   its  twisted homogeneous coordinate ring,  
 see   [Stafford \& van ~den ~Bergh 2001]  \cite{StaVdb1}.
Let $R$ be a commutative  graded ring,  such that $V=Spec~(R)$.  Denote by $R[t,t^{-1}; \sigma]$
the ring of skew Laurent polynomials defined by the commutation relation
$b^{\sigma}t=tb$
for all $b\in R$, where $b^{\sigma}$ is the image of $b$ under automorphism 
$\sigma$.  It is known, that $R[t,t^{-1}; \sigma]\cong B(V, \mathcal{L}, \sigma)$.

Let $\mathcal{H}$ be a Hilbert space and   $\mathcal{B}(\mathcal{H})$ the algebra of 
all  bounded linear  operators on  $\mathcal{H}$.
For a  ring of skew Laurent polynomials $R[t, t^{-1};  \sigma]$,  
 consider a homomorphism: 
\begin{equation}\label{eq3}
\rho: R[t, t^{-1};  \sigma]\longrightarrow \mathcal{B}(\mathcal{H}). 
\end{equation}
Recall  that  $\mathcal{B}(\mathcal{H})$ is endowed  with a $\ast$-involution;
the involution comes from the scalar product on the Hilbert space $\mathcal{H}$. 
We shall call representation (\ref{eq3})  $\ast$-coherent,   if
(i)  $\rho(t)$ and $\rho(t^{-1})$ are unitary operators,  such that
$\rho^*(t)=\rho(t^{-1})$ and 
(ii) for all $b\in R$ it holds $(\rho^*(b))^{\sigma(\rho)}=\rho^*(b^{\sigma})$, 
where $\sigma(\rho)$ is an automorphism of  $\rho(R)$  induced by $\sigma$. 
Whenever  $B=R[t, t^{-1};  \sigma]$  admits a $\ast$-coherent representation,
$\rho(B)$ is a $\ast$-algebra;  the norm closure of  $\rho(B)$  is   a   $C^*$-algebra
[Dixmier 1977]  \cite{D}.  
 We shall  denote it by ${\mathcal A}_V$ and   refer to ${\mathcal A}_V$  as   the  {\it Serre $C^*$-algebra}
 of variety $V$.

  Let $\mathcal{K}$ be the $C^*$-algebra of all compact
operators on $\mathcal{H}$.  We shall write $\tau: {\mathcal A}_V\otimes \mathcal{K}\to \mathbf{R}$
to denote   the canonical  normalized trace on  ${\mathcal A}_V\otimes \mathcal{K}$,   i.e. a positive linear functional
of norm $1$  such that $\tau(yx)=\tau(xy)$ for all $x,y\in {\mathcal A}_V\otimes \mathcal{K}$,  see  
 [Blackadar 1986] \cite[p.31]{B}.  Denote by $C(V)$ the $C^*$-algebra of complex-valued functions on the Hausdorff space $V$. 
Because ${\mathcal A}_V$  is a crossed product $C^*$-algebra of the form
${\mathcal A}_V\cong C(V)\rtimes  {\bf Z}$ \cite[Lemma 5.3.2]{N},   one can use  the Pimsner-Voiculescu 
six term exact sequence for the crossed products,  see  e.g.  [Blackadar 1986]  \cite[p.83]{B} for
 the details.  Thus   one gets the  short exact sequence of the algebraic $K$-groups:  
$0\to K_0(C(V))\buildrel  i_*\over\to  K_0({\mathcal A}_V)\to K_1(C(V))\to 0$, 
where   map  $i_*$  is induced by the natural  embedding of $C(V)$ 
into ${\mathcal A}_V$.   We  have $K_0(C(V))\cong K^0(V)$ and 
$K_1(C(V))\cong K^{-1}(V)$,  where $K^0$ and $K^{-1}$  are  the topological
$K$-groups of  $V$, see  [Blackadar 1986]  \cite[p.80]{B}. 
By  the Chern character formula,  one gets
$K^0(V)\otimes \mathbf{Q}\cong H^{even}(V; \mathbf{Q})$ and 
$K^{-1}(V)\otimes \mathbf{Q}\cong H^{odd}(V; \mathbf{Q})$, 
where $H^{even}$  ($H^{odd}$)  is the direct sum of even (odd, resp.) 
cohomology groups of $V$.  
Notice that $K_0({\mathcal A}_V\otimes \mathcal{K})\cong K_0({\mathcal A}_V)$ because
of  a stability of the $K_0$-group with respect to tensor products by the algebra 
$\mathcal{K}$,  see e.g.   [Blackadar 1986]  \cite[p.32]{B}.
One gets the   commutative diagram in Figure 1, 
where $\tau_*$ denotes  a homomorphism  induced on $K_0$ by  the canonical  trace 
$\tau$ on the $C^*$-algebra  ${\mathcal A}_V\otimes \mathcal{K}$. 
Since  $H^{even}(V):=\oplus_{i=0}^n H^{2i}(V)$ and  
$H^{odd}(V):=\oplus_{i=1}^n H^{2i-1}(V)$,   one gets  for each  $0\le i\le 2n$ 
 an injective  homomorphism   $ \tau_*:  ~H^i(V)\longrightarrow  \mathbf{R}$. 
\begin{figure}
\begin{picture}(300,100)(0,0)
\put(160,72){\vector(0,-1){35}}
\put(80,65){\vector(2,-1){45}}
\put(240,65){\vector(-2,-1){45}}
\put(10,80){$ H^{even}(V)\otimes \mathbf{Q} 
\buildrel  i_*\over\longrightarrow  K_0({\mathcal A}_V\otimes\mathcal{K})\otimes \mathbf{Q} 
\longrightarrow H^{odd}(V)\otimes \mathbf{Q}$}
\put(167,55){$\tau_*$}
\put(157,20){$\mathbf{R}$}
\end{picture}
\caption{The trace cohomology.}
\end{figure}
\begin{defn}\label{dfn2}
By an $i$-th trace cohomology  group $H^i_{tr}(V)$  of  variety  $V$   one 
understands the  abelian subgroup  of   $\mathbf{R}$ defined by the map $\tau_*$.
\end{defn}

\subsection{Langlands reciprocity}
Let $V$ be  an $n$-dimensional complex projective variety
over a number field $K$;  consider its reduction $V(\mathbf{F}_p)$ modulo
the prime ideal $\mathfrak{P}\subset K$ over a non-ramified  prime $p$.  Recall that
the {\it Weil} zeta function is defined as:  
\begin{equation}\label{eq7}
Z_p(t)=\exp~\left(\sum_{r=1}^{\infty}|V(\mathbf{F}_{p^r})|{t^r\over r}\right), \quad r\in \mathbf{C},
\end{equation}  
where $|V(\mathbf{F}_{p^r})|$ is the number  of points of variety $V(\mathbf{F}_{p^r})$
defined over the field with $p^r$ elements.  It is known that:  
\begin{equation}
Z_p(t)={P_1(t)\dots P_{2n-1}(t)\over P_0(t)\dots P_{2n}(t)},
\end{equation}  
where $P_0(t)=1-t, P_{2n}=1-p^nt$ and each $P_i(t)$ for $1\le i\le 2n-1$ is a
polynomial with integer coefficients,  such that  $P_i(t)=\prod (1-\alpha_{ij}t)$
for some algebraic integers $\alpha_{ij}$ of the absolute value $p^{i\over 2}$.  
Consider an infinite product: 
\begin{equation}\label{eq9}
L(s,V) := \prod_p  Z_p(p^{-s})={L^1(s,V)\dots L^{2n-1}(s,V)\over L^0(s,V)\dots L^{2n}(s,V)},
\end{equation}  
where $L^i(s,V)=\prod_p P_i(p^{-s})$;   the $L(s,V)$ is called 
 the {\it Hasse-Weil}   (or motivic)   $L$-function of $V$.

\medskip
On the other hand,  if $K$ is a number field  then the {\it adele ring}  ${\mathbf A}_K$
of  $K$ is a locally compact subring of the direct product $\prod K_v$ taken 
over all places $v$ of $K$;  the ${\mathbf A}_K$ is   endowed with a canonical topology.    
Consider a reductive   group $G({\mathbf A}_K)$ over 
 ${\mathbf A}_K$;   the latter  is a topological group with a canonical 
discrete subgroup $G(K)$.   
Denote by $L^2(G(K) \backslash G({\mathbf A}_K))$ the Hilbert space of all
square-integrable complex-valued functions on the homogeneous space 
$G(K)~\backslash ~G({\mathbf A}_K)$ 
and consider the right regular representation $\mathcal{R}$ of the locally compact group $G({\mathbf A}_K)$
by linear operators on the space $L^2(G(K) \backslash G({\mathbf A}_K))$
given by formula (\ref{eq1}).   It is well known, that each irreducible component
$\pi$ of the unitary representation $\mathcal{R}$ can be written in the form
$\pi=\otimes  \pi_v$, where $v$ are all unramified places of $K$.   
Using the  {\it spherical functions},  one gets an injection
$\pi_v\mapsto [A_v]$, where $[A_v]$ is a conjugacy
class of matrices in the group $GL_n(\mathbf{C})$.  The 
{\it automorphic}   $L$-function  is given by the formula:
\begin{equation}\label{eq10}
L(s,\pi)=\prod_v \left(\det ~\left[ I_n-[A_v](Nv)^{-s}\right]\right)^{-1},
\quad s\in \mathbf{C},
\end{equation}
where $Nv$ is the norm of place $v$; 
we refer the reader to [Langlands 1978] \cite[p.170]{Lan1} and   [Gelbart 1984]  \cite[p.201]{Gel1}  
for details of this construction.

\bigskip
The following conjecture relates the Hasse-Weil  and automorphic $L$-functions. 
\begin{con}\label{cnj1}
{\bf ([Langlands 1978]  \cite{Lan1})}
For each $0\le i\le 2n$ there exists an irreducible representation $\pi_i$ of the group 
$G({\mathbf A}_K)$,   such that  $L^i(s,V) \equiv L(s,\pi_i)$. 
\end{con}

\section{Proofs}
\subsection{Proof of theorem \ref{thm2}}
We shall split the proof in  two lemmas.
\begin{lem}\label{lem3.1}
The algebra ${\mathcal A}_G$ is isomorphic to a stationary AF-algebra.
 \end{lem}
\begin{proof}
Let $\mathbf{A}_K^{\times}$ be the idele group,  i.e. a group of invertible elements of the adele ring 
$\mathbf{A}_K$. Denote by $Gal~(K^{ab}|K)$ the  Galois group of the maximal abelian extension $K^{ab}$ 
of the number field $K$. The Artin reciprocity says that there exists a continuous isomorphism:
\begin{equation}\label{eq3.1}
K^{\times}\backslash\mathbf{A}_K^{\times}/C_K \longrightarrow Gal~(K^{ab}|K),
\end{equation}
where $C_K$  is the closure of the image in $K^{\times}\backslash\mathbf{A}_K^{\times}$ 
of the identity connected component of the archimedean  part $K^{\infty}$  of  the $\mathbf{A}_K^{\times}$.

Recall that $Gal~(K^{ab}|K)$ is a profinite abelian group,  i.e. a topological group isomorphic to the inverse  limit of finite abelian groups. 
It follows from the Artin reciprocity (\ref{eq3.1}),  that  $K^{\times}\backslash\mathbf{A}_K^{\times}/C_K$  is  also a profinite abelian group.  
Since every finite abelian group  is a product of the cyclic
groups ${\bf Z}/p_i^{k_i}{\bf Z}$,  we can write the group $K^{\times}\backslash\mathbf{A}_K^{\times}/C_K$ in the form:    
\begin{equation}\label{eq3.2}
K^{\times}\backslash\mathbf{A}_K^{\times}/C_K\cong\varprojlim \quad\prod_{i=1}^{l_m} \left({\bf Z}/p_i^{k_i}{\bf Z}\right), 
\end{equation}
where $m\to\infty$. Notice that the cyclic group ${\bf Z}/p_i^{k_i}{\bf Z}$ can be embedded into the finite 
field $\mathbf{F}_{q_i}$, where $q_i=p_i^{k_i}$. Thus the group $GL_n({\bf Z}/p_i^{k_i}{\bf Z})$
is correctly defined  and   from (\ref{eq3.2}) one gets an isomorphism
\begin{equation}\label{eq3.3}
GL_n(K^{\times}\backslash\mathbf{A}_K^{\times}/C_K)\cong\varprojlim \quad\prod_{i=1}^{l_m} GL_n(\mathbf{F}_{q_i}), 
\end{equation}
where $GL_n(\mathbf{F}_{q_i})$ is a finite group of order $\prod_{j=0}^{n-1}(q_i^n-q_i^j)$
and such a group  is no longer abelian. In particular, it follows from (\ref{eq3.3}) that 
the $GL_n(K^{\times}\backslash\mathbf{A}_K^{\times}/C_K)$ is a profinite group.

\bigskip
(i)  Let us show that the group $GL_n(K^{\times}\backslash\mathbf{A}_K^{\times}/C_K)$ being profinite implies that 
the  ${\mathcal A}_G$ is an AF-algebra.   Indeed,  if $G$ is 
a finite group then the  group algebra   $\mathbf{C}[G]$
has the form
\begin{equation}
\mathbf{C}[G]\cong M_{n_1}(\mathbf{C})\oplus\dots\oplus M_{n_h}(\mathbf{C}),
\end{equation}
where $n_i$ are degrees of the irreducible representations of $G$ 
and $h$  is the total number of such representations [Serre 1967]
\cite[Proposition 10]{SE}.  
In view of (\ref{eq3.3}),  we have  
\begin{equation}\label{eq3.4}
GL_n(K^{\times}\backslash\mathbf{A}_K^{\times}/C_K)\cong\varprojlim G_i,
\end{equation}
where $G_i$ is a finite group.  Consider a group algebra  
\begin{equation}\label{eq3.5}
\mathbf{C}[G_i]\cong M_{n_1}^{(i)}(\mathbf{C})\oplus\dots\oplus M_{n_h}^{(i)}(\mathbf{C})
\end{equation}
corresponding to $G_i$.  Notice that the $\mathbf{C}[G_i]$ is  a finite-dimensional 
$C^*$-algebra. The inverse limit (\ref{eq3.4})  defines an 
ascending sequence of the finite-dimensional $C^*$-algebras of the form
\begin{equation}\label{eq3.6}
\varprojlim M_{n_1}^{(i)}(\mathbf{C})\oplus\dots\oplus 
M_{n_h}^{(i)}(\mathbf{C}).  
\end{equation}
Since ${\mathcal A}_G$ is the norm closure of the group algebra 
$\mathbf{C}[GL_n (K^{\times}\backslash\mathbf{A}_K^{\times})]
\cong  \mathbf{C}[GL_n (K^{\times})\backslash GL_n(\mathbf{A}_K^{\times})]$
[Dixmier 1977]  \cite[Section 13.9]{D}, we conclude that there exists a $C^*$-homomorphism 
$h: {\mathcal A}_G\to {\Bbb A}_G$, where ${\Bbb A}_G$ is an  AF-algebra defined by 
the limit  (\ref{eq3.6}). 
To calculate the kernel of $h$, recall that $C_K\cong \varprojlim U_i$, where $U_i$ are open
subgroups of the group  $K^{\times}\backslash\mathbf{A}_K^{\times}$. We repeat the construction
of  (\ref{eq3.4})-(\ref{eq3.6}) 
and obtain an AF-algebra ${\Bbb A}_U$.  One gets an exact sequence of the
$C^*$-algebras $1\to {\Bbb A}_U\to {\mathcal A}_G\to {\Bbb A}_G\to 1$. 
In other words, the  ${\mathcal A}_G$ is an extension of the AF-algebra ${\Bbb A}_U$ 
by the AF-algebra ${\Bbb A}_G$. But any such an extension must be an AF-algebra 
itself [Brown 1982] \cite{Bro1}.  
Item (i) is proved.    

\bigskip
(ii)  It remains to prove that the ${\mathcal A}_G$ is a stationary AF-algebra.
Indeed, denote by $Fr_q$ the Frobenius map, i.e. an 
endomorphism of the finite field $\mathbf{F}_q$ acting by the formula
$x\mapsto x^q$.  The map $Fr_{q_i}$ induces  an automorphism of the group 
$GL_n(\mathbf{F}_{q_i})$.  Using  formula (\ref{eq3.3}),  one gets an automorphism 
of the group $GL_n(\mathbf{A}_K)$ and the corresponding 
group algebra $\mathbf{C}[GL_n(\mathbf{A}_K)]$. Taking the norm closure of 
the algebra $\mathbf{C}[GL_n(\mathbf{A}_K)]$,  we conclude that there exists a
non-trivial automorphism $\phi$ of the AF-algebra ${\mathcal A}_G$. 
But the AF-algebra admits   an automorphism $\phi\ne\pm ~Id$ if and only if it is  
a stationary AF-algebra [Effros 1981]  \cite[p.37]{E}.
Thus the algebra ${\mathcal A}_G$ is a stationary AF-algebra.
Lemma \ref{lem3.1} is proved.
\end{proof}

\begin{rem}\label{rmk3.3}
It follows from formula (\ref{eq3.4}) that the AF-algebra ${\mathcal A}_G$ 
is determined by a partial multiplicity matrix $B$ of rank $n$, i.e.  $B\in SL_n({\bf Z})$. 
Consider an isomorphism
\begin{equation}
{\mathcal A}_G\rtimes{\bf Z}\cong \mathcal{O}_B\otimes\mathcal{K},
\end{equation}
where the crossed product is taken by the shift automorphism of ${\mathcal A}_G$,
$\mathcal{O}_B$ is the {\it Cuntz-Krieger} algebra defined by matrix $B$ and $\mathcal{K}$
is the $C^*$-algebra of compact operators  [Blackadar 1986]  \cite[Exercise 10.11.9]{B}.
Consider a continuous group of modular  automorphisms 
$\{\sigma^t: \mathcal{O}_B\to\mathcal{O}_B ~|~t\in\mathbf{R}\}$
acting on the generators $s_1, \dots, s_n$ of the algbera  $\mathcal{O}_B$ 
by the formula $s_k\mapsto e^{it}s_k$. Then a pull back of $\sigma^t$ corresponds
to the action of  continuous symmetry group $GL_n(\mathbf{A}_K)$  on the 
homogeneous space   $GL_n(K)\backslash GL_n(\mathbf{A}_K)$.
This observation can be applied to prove Weil's conjecture on the Tamagawa numbers.  
\end{rem}
\begin{lem}\label{lem3.2}
The algebra  ${\mathcal A}_V$ embeds  into the AF-algebra ${\mathcal A}_G$,
where $V$ is a $G$-coherent variety. 
 \end{lem}
\begin{proof}
 We shall use the {\it Pimsner's Theorem} [Pimsner 1983]  \cite[Theorem 7]{Pim1}
 about an embedding of the crossed  product  algebra ${\mathcal A}_V$   
 into an  AF-algebra. It will develop that 
 the $G$-coherence of $V$ implies  that the AF-algebra  is  Morita equivalent to 
 the algebra  ${\mathcal A}_G$ of lemma \ref{lem3.1}.
We  pass to a detailed argument.

Let $V$ be a complex projective variety. 
Following    [Pimsner 1983]  \cite{Pim1} we shall think of $V$
as a compact metrizable topological space $X$. Recall that 
for a homeomorphism $\varphi: X\to X$ the point $x\in X$ 
is called {\it non-wandering} if for each neighborhood $U$ of $x$
and every $N>0$ there exists $n>N$ such that 
\begin{equation}\label{eq25}
\varphi^n(U)\cap U \ne\emptyset. 
\end{equation}
 (In other words, the point $x$ does not  ``wander''  too far from its initial 
 position in the space $X$.)   If each point $x\in X$ is a non-wandering point,
 then the homeomorphism $\varphi$ is called non-wandering.  
 
 Let $\sigma: V\to V$ be an automorphism of finite order of the $G$-coherent
 variety $V$, such that the representation (\ref{eq3}) is $\ast$-coherent. Then 
 the crossed product  
\begin{equation}\label{eq26}
{\mathcal A}_V=C(V)\rtimes_{\sigma} {\bf Z}
\end{equation}
 is the Serre $C^*$-algebra of $V$.   Since $\sigma$ is of a finite order, it is a 
 non-wandering homeomorphism of $X$.  In particular, the $\sigma$ is 
 a  pseudo-non-wandering homeomorphism   [Pimsner 1983]  \cite[Definition 2]{Pim1}.  
   Then there exists a unital (dense) embedding 
\begin{equation}\label{eq27}
{\mathcal A}_V\hookrightarrow  {\mathcal A},
\end{equation}
where ${\mathcal A}$ is an AF-algebra defined by the homeomorphism 
$\varphi$  [Pimsner 1983]  \cite[Theorem 7]{Pim1}.

Let us show that the algebra ${\mathcal A}$ is Morita equivalent
to the AF-algebra ${\mathcal A}_G$.  Indeed,  
the embedding (\ref{eq27}) induces an injective homomorphism 
of the $K_0$-groups
\begin{equation}\label{eq28}
K_0({\mathcal A}_V)\hookrightarrow  K_0({\mathcal A}).
\end{equation}
As explained in Section 2.2, the map (\ref{eq28}) defines an inclusion 
\begin{equation}\label{eq29}
H_{tr}^i(V)\subseteq K_0({\mathcal A}). 
\end{equation}
On the other hand, the trace cohomology of a $G$-coherent variety $V$
must satisfy an inclusion
\begin{equation}\label{eq30}
H_{tr}^i(V)\subseteq K_0({\mathcal A}_G). 
\end{equation}
Let $b^*=\max_{0\le i\le 2n} b_i$ be the maximal Betti number of variety $V$.
Then in formulas  (\ref{eq29}) and (\ref{eq30}) the inclusion is an isomorphism,
i.e.  $H_{tr}^*(V)\cong K_0({\mathcal A})$ and  $H_{tr}^*(V)\cong K_0({\mathcal A}_G)$.
One concludes that 
\begin{equation}\label{eq31}
 K_0({\mathcal A})\cong K_0({\mathcal A}_G).  
\end{equation}
In other words,  the AF-algebras ${\mathcal A}$ and ${\mathcal A}_G$ 
are Morita equivalent.  
The embedding ${\mathcal A}_V\hookrightarrow  {\mathcal A}_G$ follows from 
formulas (\ref{eq27}) and (\ref{eq31}).   Lemma \ref{lem3.2} is proved.
\end{proof}

Theorem \ref{thm2} follows from lemma \ref{lem3.2}.

\subsection{Proof of corollary \ref{cor1}}
Corollary \ref{cor1} follows from an observation that the Frobenius action $\sigma(Fr^i_p): H^i_{tr}(V)\to H^i_{tr}(V)$
extends to  a   Hecke operator $T_p:  K_0({\mathcal A}_G) \to K_0({\mathcal A}_G)$, 
whenever $H^i_{tr}({\mathcal A}_V)\subseteq K_0({\mathcal A}_G)$.  
Let us pass to a detailed argument.

Recall that the Frobenius map on the $i$-th trace cohomology of variety $V$
is given by an integer matrix $\sigma(Fr_p^i)\in GL_{b_i}({\bf Z})$,  where $b_i$
is the $i$-th Betti number of $V$;  moreover,
\begin{equation}\label{eq33}
|V(\mathbf{F}_p)|=\sum_{i=0}^{2n} (-1)^i ~tr~\sigma(Fr_p^i),
\end{equation}
where $V(\mathbf{F}_p)$ is the reduction of $V$ modulo a good prime $p$
\cite{Nik1}.  (Notice that (\ref{eq33})  is sufficient  to calculate the Hasse-Weil
$L$-function $L(s,V)$ of variety $V$  via equation (\ref{eq7});   hence 
the map $\sigma(Fr_p^i):   H^i_{tr}(V)\to H^i_{tr}(V)$ is motivic.) 
\begin{defn}
Denote by $T_p^i$  an endomorphism of $K_0({\mathcal A}_G)$,   such that
the diagram in Figure 2  is commutative,   where $\iota$  is the embedding (\ref{eq4}). 
By $\mathfrak{H}_i$ we understand the algebra over ${\bf Z}$  generated
by the $T_p^i\in End~ (K_0({\mathcal A}_G))$,  where $p$ runs through 
all but a finite set of primes. 
\end{defn}
\begin{figure}
\begin{picture}(300,110)(-80,0)
\put(20,70){\vector(0,-1){35}}
\put(130,70){\vector(0,-1){35}}
\put(45,23){\vector(1,0){60}}
\put(45,83){\vector(1,0){60}}
\put(0,20){$K_0({\mathcal A}_G)$}
\put(118,20){$K_0({\mathcal A}_G)$}
\put(7,80){$H^i_{tr}(V)$}
\put(115,80){$H^i_{tr}(V)$}
\put(70,30){$T_p^i$}
\put(60,90){$\sigma(Fr_p^i)$}
\put(5,50){$\iota$}
\put(140,50){$\iota$}
\end{picture}
\caption{The Hecke operator $T_p^i$.}
\end{figure}
\begin{rem}
The algebra $\mathfrak{H}_i$ is commutative.
Indeed,  the endomorphisms $T_p^i$  correspond to multiplication of the group  $K_0({\mathcal A}_G)$
by the real numbers;  the latter commute with each other. We shall call the $\{\mathfrak{H}_i ~|~ 0\le i\le 2n\}$ an $i$-th
{\it Hecke algebra}.  
\end{rem}
\begin{lem}\label{lm5}
The algebra $\mathfrak{H}_i$ defines    an  irreducible representations 
$\pi_i$ of the group $G({\mathbf A}_K)$. 
\end{lem}
\begin{proof}
Let $f\in L^2(G(K) \backslash G({\mathbf A}_K) )$
be an eigenfunction of the Hecke operators $T_p^i$;  in other words, 
the Fourier coefficients $c_p$  of the function $f$ coincide with the eigenvalues 
of the Hecke operators $T_p$ up to a scalar multiple.  Such an eigenfunction 
is defined uniquely by the algebra   $\mathfrak{H}_i$.  

Let $\mathcal{L}_f\subset  L^2(G(K) \backslash G({\mathbf A}_K) )$ 
be a subspace generated by the right  translates of $f$ by the elements of 
the locally compact group   $G({\mathbf A}_K)$.  
It is immediate (see e.g. [Gelbart 1984] \cite[Example on p. 197]{Gel1}),  
that  $\mathcal{L}_f$ is an irreducible subspace 
of the space $L^2(G(K) \backslash G({\mathbf A}_K) )$;
therefore it gives rise  to  an irreducible representation $\pi_i$ 
of the locally compact group  $G({\mathbf A}_K)$. 
Lemma \ref{lm5} follows.  
\end{proof}

\begin{lem}\label{lm6}
 $ L(s, \pi_i) \equiv L^i(s, V)$.
\end{lem}
\begin{proof}
Recall that the function $L^i(s,V)$ can be written as 
\begin{equation}\label{eq34}
L^i(s,V)=\prod_p \left(\det ~\left[ I_n-\sigma(Fr_p^i) p^{-s}\right]\right)^{-1},
\end{equation}
where $\sigma(Fr_p^i)\in GL_{b_i}({\bf Z})$ is a matrix form 
of the action of $Fr_p^i$ on the trace cohomology $H^i_{tr}(V)$.

On the other hand, from (\ref{eq10}) one gets
\begin{equation}\label{eq35}
L(s,\pi_i)=\prod_p \left(\det ~\left[ I_n-[A_p^i]p^{-s}\right]\right)^{-1},
\end{equation}  
where $[A_p^i]\subset GL_n(\mathbf{C})$ is a conjugacy class of matrices corresponding to the 
irreducible representation $\pi_i$ of the group $G({\mathbf A}_K)$. 
As explained, for such a representation we have an inclusion  
$T_p^i\in [A_p^i]$.  But  the action of the Hecke operator $T_p^i$ is an extension of 
the action of $\sigma(Fr_p^i)$ on $H_{tr}^i(V)$,   see Figure 2.   
Therefore 
\begin{equation}\label{eq36}
\sigma(Fr_p^i)=[A_p^i]
\end{equation}  
for all but a finite set of primes $p$. Comparing formulas 
(\ref{eq34})-(\ref{eq36}),  we get that   $L(s, \pi_i) \equiv L^i(s, V)$.
Lemma \ref{lm6} follows.  
\end{proof}

\bigskip
Corollary \ref{cor1} follows from lemma \ref{lm6}  and formula (\ref{eq9}).

\section{Example}
We shall illustrate theorem \ref{thm2}  and corollary \ref{cor1}
for the group $G\cong GL_2(\mathbf{A}_K)$,
where $K=\mathbf{Q}(\sqrt{D})$ is a real  quadratic field.  
\begin{prop}\label{prp1}
$K_0({\mathcal A}_G)\cong {\bf Z}+{\bf Z}\omega$,   where 
\begin{equation}\label{eq37}
\omega=
\begin{cases}
{1+\sqrt{D}\over 2},  & \hbox{if}  ~D\equiv 1 ~mod~4,\cr
               \sqrt{D},  & \hbox{if}   ~D\equiv 2,3 ~mod~4.
\end{cases}               
\end{equation}
\end{prop}
\begin{proof}
By lemma \ref{lem3.1} and remark \ref{rmk3.3}, the ${\mathcal A}_G$ 
is a stationary AF-algebra given by partial multiplicity matrix
$B\in SL_2({\bf Z})$. In particular,  $K_0({\mathcal A}_G)\cong {\bf Z}+{\bf Z}\omega$,
where $\omega\in \mathbf{Q}(\lambda_B)$, where $\lambda_B$ is the Perron-Frobenius eigenvalue
of matrix $B$. Moreover, by the construction $End~(K)\cong End~(K_0({\mathcal A}_G))$, where $End$ is the 
endomorphism ring of the corresponding ${\bf Z}$-module. But $End~(K)\cong O_K$, where $O_K$
is the ring of integers of $K$. Thus, $\lambda_B\in K$ and $\omega$ is given by formula (\ref{eq37}).  
Proposition \ref{prp1} follows.   
\end{proof}

\begin{prop}\label{prp2}
Let $\mathcal{E}_{CM}\cong \mathbf{C}/O_k$ be an elliptic curve 
with complex multiplication by the ring of integers of the imaginary quadratic 
field $k=\mathbf{Q}(\sqrt{-D})$. Then $\mathcal{E}_{CM}$ is a $G$-coherent variety of the group 
$G\cong GL_2(\mathbf{A}_K)$. 
\end{prop}
\begin{proof}
The   noncommutative torus   $\mathcal{A}_{\theta}$ is a $C^*$-algebra 
generated by the unitary operators $u$ and $v$ satisfying the commutation relation
$vu=e^{2\pi i\theta}uv$ for a constant $\theta\in\mathbf{R}$
[Rieffel 1990] \cite{Rie1}.   The Serre $C^*$-algebra of an elliptic 
curve $\mathcal{E}_{\tau}\cong  \mathbf{C}/({\bf Z}+{\bf Z}\tau)$ 
is isomorphic to $\mathcal{A}_{\theta}$ for any $\{\tau ~|~Im~\tau>0\}$, see  \cite[Theorem 1.3.1]{N}. 
In particular \cite{Nik4},  if $\tau\in O_k$ then   
\begin{equation}\label{eq38}
\begin{cases}
H_{tr}^0(\mathcal{E}_{CM}) & =H^2_{tr}(\mathcal{E}_{CM})\cong {\bf Z},\cr
           H_{tr}^1(\mathcal{E}_{CM}) & \cong   {\bf Z}+{\bf Z}\omega.
\end{cases}               
\end{equation}
Comparing formulas  (\ref{eq37}) and (\ref{eq38}),
one concludes that  $H^i_{tr}(\mathcal{E}_{CM})\subseteq K_0({\mathcal A}_G)$,
i.e. the $\mathcal{E}_{CM}$ is a $G$-coherent variety of the 
group $G\cong GL_2(\mathbf{A}_K)$. Proposition \ref{prp2}
is proved. 
\end{proof}
\begin{rem}
 The embedding of $\mathcal{A}_{\theta}$ into an AF-algebra was initially  constructed  in   
[Pimsner \& Voiculescu 1980]  \cite{PiVo1}. 
 \end{rem}
\begin{prop}\label{prp3}
$L(s, \mathcal{E}_{CM})\equiv {L(s, ~\pi_1)\over L(s, ~\pi_0)L(s, ~\pi_2)}$,
where $\pi_i$ are irreducible representations of the locally compact
 group $GL_2(\mathbf{A}_K)$.  
 \end{prop}
\begin{proof}
The Hasse-Weil $L$-function of the $\mathcal{E}_{CM}$ has the 
form: 
\begin{equation}\label{eq39}
L(s, \mathcal{E}_{CM})= {\prod_p \left[ \det~(I_2-\sigma(Fr_p^1) p^{-s}\right]^{-1}
\over \zeta(s)\zeta(s-1)},  \qquad s\in \mathbf{C}, 
\end{equation}
where $\zeta(s)$ is the Riemann zeta function and the product is taken over the set of 
good  primes;  we refer the reader to formula (\ref{eq34}).
It is immediate that  
\begin{equation}\label{eq40}
\begin{cases}
L(s, \pi_0) & =\zeta(s),\cr
L(s, \pi_2) & =\zeta(s-1),
\end{cases}               
\end{equation}
where $L(s, \pi_0)$ and $L(s, \pi_2)$ are the automorphic $L$-functions corresponding 
to the irreducible representations $\pi_0$ and $\pi_2$ of the group $GL_2(\mathbf{A}_K)$.
An irreducible representation $\pi_1$ gives rise to an automorphic $L$-function 
\begin{equation}\label{eq41}
L(s, \pi_1)= \prod_p \left(\det ~\left[ I_2-[A_p^1]p^{-s}\right]\right)^{-1}.
 \end{equation}
But formula (\ref{eq36})  says that  $[A_p^1]=\sigma(Fr_p^1)$ and therefore 
the numerator of (\ref{eq39}) coincides with the $L(s,\pi_1)$.   
Proposition \ref{prp3} is proved.
\end{proof} 
\begin{rem}
Proposition \ref{prp3}  can be proved in  terms of  the Gr\"ossencharacters 
[Silverman 1994]  \cite[Chapter II, \S 10]{S}.
\end{rem}

\subsection*{Acknowledgment}
I thank the referees for their interest and helpful comments on the draft of this paper.

\end{document}